\documentclass[11pt,reqno,letter]{amsart}

\usepackage{amssymb,amsmath,amsthm}
\usepackage{hyperref}
\usepackage{color}
\usepackage{geometry}

\theoremstyle{plain}
\newtheorem{thmM}{Theorem}

\newtheorem{thm}{Theorem}[section]

\newtheorem{lemma}[thm]{Lemma}

\theoremstyle{definition}

\newcommand{\R}{\mathbb{R}}
\newcommand{\Z}{\mathbb{Z}}
\newcommand{\N}{\mathbb{N}}

\hypersetup{
	colorlinks=true,
	citecolor=[rgb]{0,0,0.5},
	linkcolor=[rgb]{0,0,0.5},
	urlcolor=[rgb]{0,0,0.75},
	pdfpagemode=UseNone,
	pdfstartview=FitH,
	pdfdisplaydoctitle=true,
	pdftitle={Averages of simplex Hilbert transforms},
	pdfauthor={P. Durcik, J. Roos},
	pdflang=en-US
}

\DeclareFontFamily{U}{mathx}{\hyphenchar\font45}
\DeclareFontShape{U}{mathx}{m}{n}{
	<5> <6> <7> <8> <9> <10>
	<10.95> <12> <14.4> <17.28> <20.74> <24.88>
	mathx10
}{}
\DeclareSymbolFont{mathx}{U}{mathx}{m}{n}
\DeclareMathAccent{\widecheck}{0}{mathx}{"71}

\begin{document}

\title[Averages of simplex Hilbert transforms]{Averages of simplex Hilbert transforms}

\author[P. Durcik]{Polona Durcik}
\address{Polona Durcik, California Institute of Technology, 1200 E California Blvd, Pasadena CA 91125, USA}
\curraddr{Chapman University, Orange, CA, USA}
\email{durcik@caltech.edu}

\author[J. Roos]{Joris Roos}
\address{Joris Roos, University of Wisconsin-Madison, 
480 Lincoln Drive, Madison, WI 53706, USA}
\curraddr{University of Massachusetts Lowell, Lowell, MA, USA}
\email{jroos.math@gmail.com}

\date{Dec 31, 2020}

\subjclass[2010]{Primary 42B20; Secondary 42B15}  
\keywords{multilinear singular integrals, simplex Hilbert transform, paraproducts}

\numberwithin{equation}{section}
\begin{abstract}
We study a multilinear singular integral obtained by taking averages of simplex Hilbert transforms. This multilinear form is also closely related to Calder\'{o}n commutators and the twisted paraproduct. We prove $L^p$ bounds in dimensions two and three and give a conditional result valid in all dimensions. \end{abstract}

\maketitle

%%%%%%%%%%%%%%%%%%%%%%%%%%%%%%%%%%%%%%%%%%%%%%%%%%%%%%%%%%%%%%%%%%%%%%%%%%%%%%%%%%%%%%%%%%%%%%%%%%%%%%%%%%%%%%%%%

\section{Introduction}

In this paper we study the multilinear singular integral $\Lambda^{(n)}(F_0,\ldots,F_n)$ defined by
\begin{equation}\label{eqn:mainform}
\int_{\R^{n}} p.v. \int_\R \int_{[0,1]^{n-1}} F_0(x)\Big( \prod_{j=1}^{n-1} F_j(x+t\alpha_j e_j) \Big) F_n(x+te_n) d\alpha \frac{dt}t  dx,
\end{equation}
where $e_1,\ldots, e_n$ are the standard unit vectors in $\R^n$. We prove the following result.
% \begin{thmM}\label{thm:main}
% Let $p_0\in (2^{n-1},\infty]$, $p_1,\dots,p_n\in (1,\infty)$ and $\sum_{j=0}^n p_j^{-1}=1$. 

% (i) If $n=2$, then there exists a constant $c\in (0,\infty)$ such that
% \begin{equation}\label{eqn:Lpbound}
% |\Lambda^{(2)}(F_0,F_1,F_2)| \le c \prod_{j=0}^2 \|F_j\|_{p_j}.
% \end{equation}

% (ii) If $n=3$ and $p_1^{-1}+p_2^{-1}>\tfrac12$, $p_2^{-1}+p_3^{-1}>\tfrac12$, $p_1^{-1}+p_3^{-1}>\tfrac12$, then there exists a constant $c\in (0,\infty)$ such that
% \begin{equation}\label{eqn:Lpbound}
% |\Lambda^{(3)}(F_0,F_1,F_2,F_3)| \le c \prod_{j=0}^3 \|F_j\|_{p_j}.
% \end{equation}
% The constant $c$ depends on $p_0,p_1,\dots,p_n$.
% \end{thmM}
\begin{thmM}\label{thm:main}
Let $n$ equal two or three. 
Let $p_0\in (2^{n-1},\infty)$, $p_1,\dots,p_n\in (1,\infty)$,  and $\sum_{j=0}^n p_j^{-1}=1$. If $n=3$, assume in addition that $\max(p_1,p_2,p_3)<4$. Then there exists a constant $c\in (0,\infty)$ such that
\begin{equation}\label{eqn:Lpbound}
|\Lambda^{(n)}(F_0,\dots,F_n)| \le c \prod_{j=0}^n \|F_j\|_{p_j}.
\end{equation}
The constant $c$ depends on $p_0,p_1,\dots,p_n$. The estimate fails if $p_j=\infty$ for some $j~=~1,\dots,n$.
\end{thmM}
If we drop the average over $\alpha$ in \eqref{eqn:mainform} and fix, say $\alpha=(1,\dots,1)$, then we obtain the \emph{simplex Hilbert transform}
\begin{equation}\label{eqn:simplexht}
\int_{\R^{n}} p.v. \int_\R  F_0(x)\Big( \prod_{j=1}^{n-1} F_j(x+t e_j) \Big) F_n(x+te_n)  \frac{dt}t  dx.
\end{equation}
Note that for $n=1$, both \eqref{eqn:mainform} and \eqref{eqn:simplexht} reduce to a dualized form of the classical Hilbert transform. For $n\ge 2$, it is currently not known whether the multilinear form \eqref{eqn:simplexht} satisfies any $L^p$ estimates. This is a major open problem in harmonic analysis. In the case $n=2$, the form \eqref{eqn:simplexht} is also called \emph{triangular Hilbert transform} and some partial progress has been made in \cite{KTZ15}, where bounds for a Walsh model are obtained, when one of the three functions takes a special form. $L^p$ bounds for the triangular Hilbert transform would unify several known results in time-frequency analysis, such as uniform bounds for the bilinear Hilbert transform and bounds for Carleson's operator, which controls pointwise almost-everywhere convergence of Fourier series. More generally, bounds for the $n$-simplex Hilbert transform would imply bounds for the \emph{$n$-linear Hilbert transform}
\begin{equation}\label{eqn:multilinht}
\int_{\R} p.v. \int_{\R} f_0(x) f_1(x+t) f_2(x+2t)\cdots f_n(x+nt) \frac{dt}t  dx.
\end{equation}
We refer to \cite[Appendix B]{KTZ15} and \cite[Appendix A]{Zor17} for the details of this implication. If $n=2$, then the previous display represents a trilinear form dual to the bilinear Hilbert transform, which was first shown to satisfy $L^p$ bounds by Lacey and Thiele \cite{LT97}, \cite{LT99}. If $n\ge 3$, then it is not known if \eqref{eqn:multilinht} satisfies any $L^p$ bounds and this is considered to be a difficult open problem. An averaged form of the $3$-linear Hilbert transform was studied by Palsson \cite{Pal12}. Tao \cite{Tao16} showed that truncations of the multilinear form \eqref{eqn:multilinht} exhibit some non-trivial cancellation by using regularity lemmas from arithmetic combinatorics. For the more general simplex Hilbert transform \eqref{eqn:simplexht}, a similar result was obtained by Zorin-Kranich \cite{Zor17} and this was strengthened subsequently in \cite{DKT16}.

The multilinear form \eqref{eqn:mainform} is closely related to the  forms dual to \emph{Calder\'{o}n commutators}, which can be written as 
\begin{equation}\label{eqn:commutators}
\int_{\R}  p.v. \int_\R  \int_{[0,1]^{n-1}} f_0(x)\Big( \prod_{j=1}^{n-1} f_j(x+t\alpha_j) \Big) f_n(x+t) d\alpha  \frac{dt}t  dx
\end{equation}
for a suitable choice of the functions $f_j$. 
In fact, $L^p$ bounds for \eqref{eqn:mainform} imply the same $L^p$ bounds for \eqref{eqn:commutators}. This follows in the same way as bounds for the simplex Hilbert transform would imply bounds for \eqref{eqn:multilinht}. 
%. This can be seen by choosing the functions $F_j$ to be appropriate elementary tensors (see \cite[Appendix A]{Zor17}).
Calder\'{o}n commutators played a key role in the resolution of the problem of bounding the Cauchy integral along Lipschitz curves by Coifman, McIntosh and Meyer \cite{CMM82}. More recently, Muscalu \cite{Mus14a}, \cite{Mus14b} developed an alternative approach for bounding Calder\'{o}n commutators, which also inspired our proof of Theorem \ref{thm:main}. We remark that in view of the connection with the Cauchy integral, the focus for bounding commutators often lies in obtaining $L^\infty$ bounds and a constant that depends polynomially on $n$. In sharp contrast to this, the estimate \eqref{eqn:Lpbound} fails if $p_j=\infty$ for some $j\in \{1,\dots,n\}$. We prove this in \S \ref{sec:counterex}.

Another object that is intimately related to \eqref{eqn:mainform} is the \emph{twisted paraproduct}
\begin{equation}\label{eqn:twisted}
\int_{\R^4} F_0(x,y) F_1(x-t,y) F_2(x,y-s) K(t,s) d(t,s,x,y),
\end{equation}
where $K$ is a Calder\'{o}n--Zygmund kernel on $\R^2$. This operator first arose as a degenerate case of the two-dimensional bilinear Hilbert transform that did not fall within the scope of Demeter and Thiele's time-frequency techniques \cite{DT10}. $L^p$ bounds were obtained by Kova\v{c} \cite{K12:tp}, who used a novel argument involving repeated use of the Cauchy--Schwarz inequality and a telescoping identity. These techniques have since been further developed and applied to other problems (including the application to the simplex Hilbert transform \cite{DKT16} which we have already mentioned).

Let $m$ be a bounded function on $\R^n\setminus\{0\}$. Then we define
\begin{equation}\label{eqn:lambdam}
\Lambda_m (F_0,\ldots F_n) = \int_{(\R^n)^n} \widehat{F_0}(-(\tau_1+\cdots + \tau_n))\Big(\prod_{j=1}^n \widehat{F_j}(\tau_j) \Big) m(\tau_{11},\dots,\tau_{nn})  d\tau,
\end{equation}
where $\tau=(\tau_1,\dots,\tau_n)\in(\R^{n})^n$ and $\tau_j = (\tau_{j1},\dots, \tau_{jn})\in \R^n$. The key property here is that the multiplier only depends on the diagonal frequency variables $\tau_{11},\dots,\tau_{nn}$. The form \eqref{eqn:twisted} is obtained by specializing $n=2$.  If we set
\begin{equation}\label{eqn:multiplier}
\mu_n(\xi) = \int_{[0,1]^{n-1}} \mathrm{sgn}(\alpha_1\xi_1 + \cdots + \alpha_{n-1} \xi_{n-1} + \xi_n) d\alpha,
\end{equation}
then $\Lambda_{\mu_n}$ conincides with our multilinear form \eqref{eqn:mainform} (up to a multiplicative constant). The symbol \eqref{eqn:multiplier} also coincides with the symbol of the Calder\'{o}n commutator \eqref{eqn:commutators}. The function $\mu_n$ is continuous on $\R^n\setminus\{0\}$ but fails to be differentiable on several hyperplanes. Some details and further comments on the function $\mu_n$ are contained in \S \ref{sec:multiplier}.

We call a bounded function $m\in C^{\infty}(\R^n\setminus\{0\})$ a \emph{standard smooth symbol} if for every $\alpha\in\N^n_0$ there exists $c_\alpha\in (0,\infty)$ such that for every $\xi\in\R^n\setminus\{0\}$ we have
\begin{equation}\label{eqn:hmest}
|\partial^{\alpha} m(\xi)| \le c_\alpha |\xi|^{-|\alpha|}.
\end{equation}
Notice that $\mu_n$ is not a standard smooth symbol. 
Moreover, note that even if $m$ is a standard smooth symbol, the multilinear Fourier multipliers of the form \eqref{eqn:lambdam} do not belong to the multilinear Calder\'{o}n--Zygmund paradigm as set forth in \cite{GT02}. 

Our proof of Theorem \ref{thm:main} relies on the estimate
\begin{equation}\label{eqn:Lpboundhm}
 |\Lambda_m(F_0,\ldots,F_n)| \le c \prod_{j=0}^n \|F_j\|_{p_j}
 \end{equation}
for standard smooth symbols $m$ with the constant $c$ only depending on $m$. In fact, we prove the following conditional result.

\begin{thmM}\label{thm:cond}
Let $n\ge 2$. If the estimate \eqref{eqn:Lpboundhm} holds for all standard smooth symbols $m$ and some fixed exponents $p_0,\dots,p_n\in (1,\infty)$, then \eqref{eqn:Lpbound} holds (for the same exponents).
\end{thmM}

In the case $n=2$, Kova\v{c} \cite{K12:tp} showed that \eqref{eqn:Lpboundhm} holds for all $p_0\in (2,\infty), p_1, p_2\in (1,\infty)$ satisfying $p_0^{-1}+p_1^{-1}+p_2^{-1}=1$.  

In the case $n=3$ we are still able to use the techniques from \cite{K12:tp} to prove \eqref{eqn:Lpboundhm} for all $p_0\in (4,\infty), p_1,p_2,p_3\in (1,4)$ satisfying $p_0^{-1}+p_1^{-1}+p_2^{-1}+p_3^{-1}=1$. Indeed, the proof is by reduction to a result for a certain dyadic form considered in \cite{K12:tp}, combined with a transition from the dyadic to the continuous case and fiber-wise Calder{\'o}n--Zygmund decomposition by Bernicot \cite{Ber12}. These arguments are well-known, but for convenience of the reader some details are given in \S \ref{sec:twisted}.

We currently do not have the estimate \eqref{eqn:Lpboundhm} for any $n\geq 4$ (for any tuple of exponents $(p_j)_j$). This problem is closely related to the question of extending the range of exponents for the twisted paraproduct \eqref{eqn:twisted}, as well as proving $L^p$ estimates  for the integrals studied in \cite{DT18}, in the case when the index set is an arbitrary subset of the cube in $\R^n$. However, in \S \ref{sec:main} we show that \eqref{eqn:Lpbound} holds provided that \eqref{eqn:Lpboundhm} holds.

{\noindent\em Organization of the paper.} In \S \ref{sec:main} we prove Theorem \ref{thm:cond}. In \S \ref{sec:twisted} we prove the estimate \eqref{eqn:Lpboundhm} for $n=3$. In \S \ref{sec:counterex} we give a counterexample to \eqref{eqn:Lpbound} if $p_j=\infty$ for some $j\in\{1,\dots,n\}$. In \S \ref{sec:multiplier} we make some comments about the symbol $\mu_n$.\\

{\noindent\em Acknowledgment.} We thank Vjeko Kova\v{c} for helpful comments on a preliminary draft and Lenka Slav\'{i}kov\'{a} for pointing out an oversight in a previous version of this paper.

\section{Main argument}\label{sec:main}
In this section we fix exponents $p_0,\dots,p_n\in (1,\infty)$ and prove that \eqref{eqn:Lpbound} holds under the assumption that for all standard smooth symbols $m$ we have
\begin{equation}\label{eqn:hmassumption}
 |\Lambda_m(F_0,\ldots,F_n)| \lesssim \prod_{j=0}^n \|F_j\|_{p_j},
\end{equation}
where $\Lambda_m$ is defined in \eqref{eqn:lambdam}.
%If $n=2$, the estimate \eqref{eqn:hmassumption} is proven in \cite{K12:tp}, while the case $n=3$ is discussed in \S \ref{sec:twisted} below. If $n\ge 4$, we do not have \eqref{eqn:hmassumption}, but we can still prove \eqref{eqn:Lpbound} under this assumption.

First we make some comments about terminology. We say that a function $\psi$ of one real variable is \emph{of $\psi$-type} if it is a smooth function compactly supported in an interval not containing the origin. The letter $\psi$ will always denote a $\psi$-type function supported on $[-2,-\frac12]\cup [\frac12,2]$ such that $\sum_{k\in\Z} \psi(2^{-k} \eta)=1$ for all $\eta\not=0$. We also reserve the letters $\widetilde{\psi}, \psi^1, \psi^2, \dots$ for other $\psi$-type functions, which may change throughout the text.

Similarly, a function $\varphi$ of one real variable is \emph{of $\varphi$-type} if it is a smooth function that is compactly supported. The letters $\varphi, \widetilde{\varphi}, \varphi^{1}, \varphi^{2}, \dots$ are reserved for $\varphi$-type functions, which may change at different stages of the argument.

If $\varphi$ is a function and $k\in\Z$ an integer, then we write $\varphi_k(\eta)=\varphi(2^{-k}\eta)$.

We start by decomposing $\R^n$ into a family of Whitney boxes. Then we have
\[ 1 = \sum_{k_1,\dots,k_n\in\Z} \psi_{k_1}(\xi_1) \cdots \psi_{k_n}(\xi_n) \]
for almost every $\xi\in\R^n$. We split the sum over $(k_1,\dots,k_n)\in\Z^n$ according to which of the components $k_j$ is the largest:
\begin{equation}\label{eqn:partunity0}
1 = \sum_{k_1\ge k_2,\dots,k_n} (\cdots) + \sum_{k_2\ge k_1,k_3,\dots,k_n, k_2\not=k_1} (\cdots) + \cdots + \sum_{k_n>k_1,\dots,k_{n-1}} (\cdots).
\end{equation}
For each of the $n$ terms in the previous display, we perform a further decomposition. Here we have that $k_1$ is no smaller than $k_2,\dots,k_n$. We want to distinguish the case that $k_1$ is \emph{much} larger than all of the remaining integers $k_2,\dots,k_n$ from the case that $k_1$ is roughly equal to at least one of the numbers $k_2,\dots,k_n$. To be more precise, we fix a large integer $\ell_0\gg 1$. For example, $\ell_0=100\log(n)$ will do. Then decompose the corresponding summation as
\begin{equation}\label{eqn:partunity1}
\sum_{k_1\ge k_2,\dots,k_n} (\cdots) = \sum_{k_1-\ell_0\ge k_2,\dots,k_n} (\cdots)  + \sum_{\substack{k_1\ge k_2,\dots,k_n,\\\exists j\in\{2,\dots,n\}\,\mathrm{s.t.}\,k_1-\ell_0 < k_j }} (\cdots).
\end{equation}
We say that the first of these two terms is \emph{of $\psi$-$\varphi$-type}. This is because we can sum in $k_2,\dots,k_n$ to obtain
\[ \sum_{k_1-\ell_0\ge k_2,\dots,k_n} \psi_{k_1}(\xi_1) \cdots \psi_{k_n}(\xi_n) = \sum_{k_1\in\Z} \psi_{k_1}(\xi_1) \prod_{j=2}^n \sum_{k_j\le k_1+\ell_0} \psi_{k_j}(\xi_j)\] 
\[= \sum_{k\in\Z} \psi_k(\xi_1) \varphi_{k}(\xi_2)\cdots \varphi_{k}(\xi_n). \]
Similarly, the second term in \eqref{eqn:partunity1} can be written as a sum over $O(n)$ terms of the form
\[ \sum_{k\in\Z} \psi^1_k(\xi_1) \psi^{j_1}_{k}(\xi_{j_1}) \prod_{j\not=j_1} \varphi^{j}_k(\xi_j), \]
where $\psi^j, \varphi^j$ are $\psi$- and $\varphi$-type functions, respectively. We refer to such a term as being \emph{of $\psi$-$\psi$-type}.

Summarizing, we obtain a partition of unity into $O(n^2)$ terms each of which is either of $\psi$-$\varphi$-type or of $\psi$-$\psi$-type.
Now we write our multilinear form as
\begin{equation} \label{eqn:offcoordpf0} \int_{(\R^n)^n} \widehat{F_0}(-(\tau_1+\cdots+\tau_n))\widehat{F_1}(\tau_1) \cdots \widehat{F_n}(\tau_n) \mu_n(\tau_{11},\dots,\tau_{nn}) d(\tau_1,\dots,\tau_n),
\end{equation}
and split the symbol $\mu_n$ into $O(n^2)$ pieces according to the aforementioned partition of unity. At this point our argument splits into three cases.

{\em Case 1: the $\psi$-$\psi$ case}. Symbols of $\psi$-$\psi$-type can be handled using an argument similar to Muscalu's alternative proof \cite{Mus14b} of the Coifman--McIntosh--Meyer theorem. This is detailed in \S \ref{sec:psipsi}.

{\em Case 2: the smooth $\psi$-$\varphi$ case.} For the $\psi$-$\varphi$-type pieces we distinguish two separate cases. The first scenario is that the $\psi$ falls onto the $n$th frequency variable (corresponding to the last term in \eqref{eqn:partunity0}). Here we are concerned with the symbol
\begin{equation}\label{eqn:psiphismooth}
\sum_{k\in\Z} \mu_n(\xi)\varphi_k(\xi_1) \cdots \varphi_{k}(\xi_{n-1}) \psi_{k}(\xi_n),
\end{equation}
where $\varphi$ is supported in $[-2^{-\ell_0+1}, 2^{-\ell_0+1}]$. Recalling that $\psi$ is supported in $[-2,-\tfrac12]\cup [\tfrac12,2]$ and $\ell_0$ is large, observe that the support of \eqref{eqn:psiphismooth} is disjoint from each of the bad hyperplanes
\[ \{\xi\in\R^n\,:\, \alpha_1 \xi_1 + \cdots + \alpha_{n-1}\xi_{n-1} + \xi_n=0 \}\quad (\alpha\in \{0,1\}^{n-1}). \]
Thus, $\mu_n$ is smooth on the support of \eqref{eqn:psiphismooth} (see Lemma \ref{lem:hyperplanes}) and therefore the desired bounds for the multilinear form with symbol \eqref{eqn:psiphismooth} follow from \eqref{eqn:hmassumption}.

{\em Case 3: the rough $\psi$-$\varphi$ case.} It remains to consider $\psi$-$\varphi$-type pieces where $\psi$ falls onto one of $\xi_1$, \dots, $\xi_{n-1}$. By symmetry, we may restrict our attention to the symbol
\begin{equation}\label{eqn:psiphin}
\sum_{k\in\Z} \mu_n(\xi)\psi_k(\xi_1) \varphi_k(\xi_2) \cdots \varphi_{k}(\xi_{n}),
\end{equation}
where $\varphi$ is supported in a small neighborhood of the origin. The support of this symbol is contained in a narrow cone around the $\xi_1$-axis and it has significant intersection with the singularities of $\mu_n$. In this case the idea is to reduce to a symbol that vanishes on the $\xi_1$-axis and then consider an appropriate lacunary decomposition with respect to that axis, where each of the pieces is a $\psi$-$\psi$-type symbol which can be handled as in Case 1. This argument is detailed in \S \ref{sec:psiphirough}. We note that by inspection of the multiplier a simpler argument is available if $n=2$. 

\subsection{\texorpdfstring{The $\psi$-$\psi$ case}{The Psi-Psi case}}\label{sec:psipsi}

%\begin{lem}
%Let $1\le j_0\not=j_1\le n$. Consider the symbol
%\[ m(\xi) = \sum_{k\in\Z} \mu_n(\xi) \psi^{j_0}_k(\xi_{j_0}) \psi^{j_1}_k(\xi_{j_1})\prod_{j\not=j_0,j_1} \varphi^{j}_k(\xi_j) \]
%\end{lem}

In this section we bound the multilinear form corresponding to a symbol of the form
\[ \sum_{k\in\Z} \mu_n(\xi) \psi^{j_0}_k(\xi_{j_0}) \psi^{j_1}_k(\xi_{j_1})\prod_{j\not=j_0,j_1} \varphi^{j}_k(\xi_j), \]
where $1\le j_0<j_1\le n$ and the $\psi^j, \varphi^j$ are of $\psi$- and $\varphi$-type, respectively.\\
Based on the symmetries of the function $\mu_n$, we need to distinguish two cases: $(j_0,j_1)=(1,2)$ and $(j_0,j_1)=(1,n)$ (they are different only if $n>2$). However, the arguments we use for these two cases are, \emph{mutatis mutandis}, identical. Thus, we will restrict our attention here to the case $(j_0,j_1)=(1,n)$. To improve readability we also suppose that the various $\psi$- and $\varphi$-type functions are identical. This is again no loss of generality. With these reductions in mind, we are now left with the symbol
\[ \sum_{k\in\Z} \mu_n(\xi) \psi_k(\xi_1) \varphi_k(\xi_2)\cdots \varphi_k(\xi_{n-1}) \psi_k(\xi_n),\]
where $\psi$ is a smooth function supported in $[1/2, 2]\cup [-2,-1/2]$ and $\varphi$ is a smooth function supported in $[-2,2]$. Recall that
\[ \mu_n(\xi) = \int_{[0,1]^{n-1}} \mathrm{sgn}(\alpha_1 \xi_1+\cdots+\alpha_{n-1} \xi_{n-1} + \xi_n) d\alpha \]

Our proof follows the argument of Muscalu \cite{Mus14b}. The argument simplifies because we do not need to make use of the non-compact Littlewood--Paley projections used in \cite{Mus14b}. The key observation that we will use (due to \cite{Mus14b}) is that
\[ \mu_n(\xi) = \int_{[0,1]^{n-2}} \mu_2(\xi_1, \alpha_2 \xi_2 + \cdots + \alpha_{n-1}\xi_{n-1}+\xi_n) d(\alpha_2,\dots,\alpha_{n-1}).  \]
This will allow us to make use of the fact that the Fourier coefficients of $\mu_2$, when restricted to a Whitney box, decay quadratically ( see Lemma \ref{lem:fourierdecay} below). For $\mu_n$ we seem to get only linear decay (see \cite{Mus14a}). We introduce the variable
\[ \eta = \eta_{\alpha,\xi} = \alpha_2 \xi_2 + \cdots + \alpha_{n-1} \xi_{n-1} + \xi_n \]
and rewrite the quantity $\mu(\xi)$ as
\begin{equation}\label{eqn:offcoordpf1} \int_{[0,1]^{n-2}}  \sum_{k\in\Z} \rho(2^{-k}(\xi_1,\eta)) \psi_k(\xi_1) \varphi_k( \xi_2)\cdots \varphi_k(\xi_{n-1}) \psi_k(\xi_n)d(\alpha_2,\dots,\alpha_{n-1}),
\end{equation}
where we have set
\[ \rho(\xi_1,\eta) = \mu_2(\xi_1,\eta) \widetilde{\psi}(\xi_1)\widetilde{\varphi}(\eta).\]
By the Fourier inversion formula we have
\begin{equation}\label{eqn:offcoordpf2}
\rho(\xi_1, \eta) = \int e^{-2\pi i (\xi_1 u + \eta v)} \widecheck{\rho}(u,v) d(u,v).
\end{equation}

At this point we record the following estimate for the function $\widecheck{\rho}$.
\begin{lemma}\label{lem:fourierdecay}
For every $N\in\N, u,v\in\R$ we have
\begin{equation}\label{eqn:quadfourierdecay}
|\widecheck{\rho}(u,v)| \lesssim_N (1+|v|)^{-2} ((1+|v-u|)^{-N} + (1+|u|)^{-N}).
\end{equation}
\end{lemma}
This is obtained by writing $\widecheck{\rho}(u,v)=\int e^{2\pi i (\xi_1 u + \eta v)} \mu_2(\xi_1,\eta)\widetilde{\psi}(\xi_1)\widetilde{\varphi}(\eta)d(\xi_1,\eta)$ and first integrating by parts twice with respect to $\eta$ and then $N$ times with respect to $\xi_1$. For the details we refer to \cite[Lemma 2.4]{Mus14a}. Alternatively, the reader may look ahead at the details of our proof of \eqref{eqn:kernelest} below, which is a slight variant of this estimate.

Note that the right hand side of \eqref{eqn:quadfourierdecay} is integrable in $u,v\in\R$. Let us for the moment fix $\alpha_2,\dots,\alpha_{n-1}\in [0,1]$, $u,v\in\R$ and $k\in\Z$. After expanding $\rho(\xi_1,\eta)$ according to \eqref{eqn:offcoordpf2}, the corresponding contribution to our multilinear form \eqref{eqn:offcoordpf0} becomes
\[\nonumber \int_{(\R^n)^n} \Big(\widehat{F_1}(\tau_1)\psi_k(\tau_{11}) e^{-2\pi i 2^{-k} \tau_{11} u}\Big) \prod_{j=2}^{n-1}\Big(\widehat{F_j}(\tau_j)\varphi_k(\tau_{jj}) e^{-2\pi i 2^{-k} \alpha_{jj}\tau_{jj} v}\Big) \]
\begin{equation}\label{eqn:offcoordpf3}
\Big(\widehat{F_n}(\tau_n)\psi_k(\tau_{nn})e^{-2\pi i 2^{-k} \tau_{nn} v} \Big) \widehat{F_0}(-(\tau_1+\cdots+\tau_n))d(\tau_1,\dots,\tau_n), \end{equation}
Expanding the Fourier transform $\widehat{F_0}$, the previous display becomes
\[\nonumber \int_{\R^n} \Big(\int_{\R^n}\widehat{F_1}(\tau_1)\psi_k(\tau_{11}) e^{-2\pi i (2^{-k} \tau_{11} u ) } e^{2\pi i x\cdot \tau_1} d\tau_1\Big) \prod_{j=2}^{n-1}\Big( \int_{\R^n} \widehat{F_j}(\tau_j)\varphi_k(\tau_{jj}) e^{-2\pi i 2^{-k} \alpha_{jj} \tau_{jj} v} e^{2\pi i x\cdot \tau_j} d\tau_j \Big) \]
\[ \Big(\int_{\R^n} \widehat{F_n}(\tau_n)\psi_k(\tau_{nn})e^{-2\pi i 2^{-k} \tau_{nn} v} e^{2\pi i x\cdot \tau_{n}} d\tau_n \Big)
F_0(x)dx, \]
which by the Fourier inversion formula is equal to
\[ \int_{\R^n} 
F_0(x)  (F_1*_1 \widecheck{\psi_k}^{(u)})(x)  \prod_{j=2}^{n-1} (F_j*_j \widecheck{\varphi_k}^{(\alpha_j v)})(x) (F_n *_n \widecheck{\psi_k}^{(v)})(x) dx. \]
Here $*_j$ denotes convolution in the $j$th variable and $f^{(u)}(x) = f(x+u)$ denotes translation.
Summing in $k\in\Z$ we obtain the multilinear form
\begin{equation}\label{eqn:offcoordpf4} \int_{\R^n} F_0(x)\sum_{k\in\Z} 
(F_1*_1 \widecheck{\psi_k}^{(u)})(x)  \prod_{j=2}^{n-1} (F_j*_j \widecheck{\varphi_k}^{(\alpha_j v)})(x) (F_n *_n \widecheck{\psi_k}^{(v)})(x) dx.
\end{equation}
Using the triangle inequality and the Cauchy--Schwarz inequality we estimate the absolute value of \eqref{eqn:offcoordpf4} by
\[\nonumber \int_{\R^n} |F_0(x)| \Big(\sum_{k\in\Z} |(F_1*_1 \widecheck{\psi_k}^{(u)})(x)|^2\Big)^{1/2} \prod_{j=2}^{n-1} \sup_{k\in\Z} |(F_j*_j \widecheck{\varphi_k}^{(\alpha_j v)})(x)|\]
\begin{equation}\label{eqn:offcoordpf6}\Big(\sum_{k\in\Z} |(F_n *_n \widecheck{\psi_k}^{(v)})(x)|^2\Big)^{1/2}   dx \end{equation}

At this point we need to make use of the following well-known bounds for \emph{shifted maximal functions} and \emph{shifted square-functions}.
\begin{lemma}
Let $u\in\R$, $\varphi$ a smooth bump function and $\psi$ a $\psi$-type function. Then we have
\begin{equation}\label{eqn:shiftedmax}
\| \sup_{k\in\Z} |f*\widecheck{\varphi_k}^{(u)}| \|_p \lesssim \log(2+|u|) \|f\|_p 
\end{equation}
and
\begin{equation}\label{eqn:shiftedsqfct}
\Big\| \Big(\sum_{k\in\Z} |f*\widecheck{\psi_k}^{(u)}|^2 \Big)^{1/2} \Big\|_p \lesssim \log(2+|u|) \|f\|_p
\end{equation}
for every $p\in (1,\infty)$.
\end{lemma}
For the proof of \eqref{eqn:shiftedmax} we refer to \cite[Theorem 4.1]{Mus14a} (it is also contained in \cite{Ste93}) and for the proof of \eqref{eqn:shiftedsqfct} see \cite[Theorem 5.1]{Mus14a}.
Using H\"older's inequality, \eqref{eqn:shiftedmax}, \eqref{eqn:shiftedsqfct} and Fubini's theorem, we see that \eqref{eqn:offcoordpf6} is bounded by a constant times
\begin{equation}\label{eqn:offcoordpf5}
\prod_{j=0}^n \|F_j\|_{p_j}
\end{equation}
for all tuples $(p_0,\dots,p_n)$ satisfying $\sum_{j=0}^n p_j^{-1}=1$, $p_j\in (1,\infty)$ ($p_j=\infty$ is also allowed unless $j\in\{1,n\}$) and the constant takes the form
\[ C \log(2+|u|) \log(2+|v|) , \]
where $C$ is an absolute constant. Integrating over $u,v\in\R$, $\alpha_2,\dots,\alpha_{n-1}\in [0,1]$ and making use of Lemma \ref{lem:fourierdecay} we obtain the desired bounds for our multilinear form.

\subsection{\texorpdfstring{The rough $\psi$-$\varphi$ case}{The rough Psi-Phi case}}\label{sec:psiphirough}
In this section we deal with the symbol
\[\sum_{k\in\Z} \mu_n(\xi)\psi_k(\xi_1) \varphi_k(\xi_2) \cdots \varphi_{k}(\xi_{n}).\]
The difficulty is that the support of this function intersects some of the hyperplanes on which $\mu_n$ fails to be differentiable.  The first step is again to use
\[ \mu_n(\xi) = \int_{[0,1]^{n-2}} \mu_2(\xi_1, \alpha_2 \xi_2 + \cdots + \alpha_{n-1}\xi_{n-1}+\xi_n) d(\alpha_2,\dots,\alpha_{n-1}).  \]
We again consider the variable
\[ \eta = \eta_{\alpha,\xi} = \alpha_2 \xi_2 + \cdots + \alpha_{n-1} \xi_{n-1} + \xi_n. \]
Fixing $\alpha_2,\dots,\alpha_{n-1}\in [0,1]$, we now study the symbol
\[ \sum_{k\in\Z} \mu_2(\xi_1, \eta) \psi_k(\xi_1) \prod_{j=2}^n \varphi_k(\xi_j).\]

We split $\psi=\psi^+ + \psi^-$ such that $\psi^+$ is supported on $\{\eta>0\}$ and $\psi^-$ is supported on $\{\eta<0\}$. Note that $\psi^+$ and $\psi^-$ are smooth. We split the symbol accordingly into two pieces. By symmetry it suffices to consider the symbol
\[ \sum_{k\in\Z} \mu_2(\xi_1, \eta) \psi^+_k(\xi_1) \prod_{j=2}^n \varphi_k(\xi_j).\]
Using \eqref{eqn:hmassumption}, we may subtract the smooth symbol $\sum_{k\in\Z} \psi^+_k(\xi_1) \prod_{j=2}^n \varphi_k(\xi_j)$ so that we are now concerned with
\[ \sum_{k\in\Z} \widetilde{\mu}_2(\xi_1,\eta) \psi^+_k(\xi_1) \prod_{j=2}^n \varphi_k(\xi_j),\]
where we have set
\[ \widetilde{\mu}(\xi_1,\eta)= \mu_2(\xi_1, \eta)-1  = -2\int_0^1 \mathbf{1}_{<0}(\xi_1\alpha+\eta) d\alpha.\]

Now we perform another paraproduct decomposition. More precisely, we expand each of the $\varphi_k(\xi_j)$ into an appropriate sum of $\psi$-type functions and distinguish $n-1$ terms based on which of the summation indices is the largest (as in \eqref{eqn:partunity0}). With this in mind it suffices to consider the symbol
\[ \sum_{\ell>\ell_0} \Big( \sum_{k\in\Z} \widetilde{\mu}(\xi_1,\eta) \psi^+_k(\xi_1) \psi_{k-\ell}(\xi_2) \prod_{j=3}^{n} \varphi_{k-\ell}(\xi_j)  \Big) \]

In the following we will fix $\ell>\ell_0$ and bound the multilinear form corresponding to the inner sum with summable decay in $\ell$. Shifting the index $k$ and using homogeneity of $\mu$ we can write the inner sum as
\begin{equation}\label{eqn:psiphiroughpf1} \sum_{k\in\Z} \widetilde{\mu}(2^{-k}(\xi_1,\eta)) \psi^+_{k+\ell}(\xi_1) \psi_{k}(\xi_2) \prod_{j=3}^{n} \varphi_{k}(\xi_j)
\end{equation}
Observe that on the support of each summand we have 
\[|\eta| \le |\xi_2| + \cdots + |\xi_{n-1}|+|\xi_n|\le 2n\cdot 2^{-k}. \]
Choose a smooth function $\widetilde{\varphi}$ such that it equals $1$ on the interval $[-2n,2n]$ and is equal to $0$ on a slightly larger interval. Similarly, we choose a smooth function $\widetilde{\psi}$ which is $1$ on the interval $[\frac12,2]$ and equal to $0$ on, say, $[\frac12-\frac1{100},2+\frac1{100}]$. Then \eqref{eqn:psiphiroughpf1} is equal to
\[ \sum_{k\in\Z} \widetilde{\mu}(2^{-k}(\xi_1,\eta))\widetilde{\psi}_{k+\ell}(\xi_1) \widetilde{\varphi}_k(\eta) \psi^+_{k+\ell}(\xi_1) \psi_{k}(\xi_2) \prod_{j=3}^{n} \varphi_{k}(\xi_j),\]
which we rewrite as
\begin{equation}\label{eqn:psiphiroughpf3}
\sum_{k\in\Z} m_\ell(2^{-k-\ell}\xi_1, 2^{-k}\eta) \psi_{k+\ell}(\xi_1) \psi_{k}(\xi_2) \prod_{j=3}^{n} \varphi_{k}(\xi_j),
\end{equation}
where
\[ m_\ell(\xi_1,\eta) = \widetilde{\mu}(\xi_1,2^{-\ell}\eta)\widetilde{\psi}(\xi_1) \widetilde{\varphi}(\eta). \]
Here we have used homogeneity of $\widetilde{\mu}$. By the Fourier inversion formula we have
\begin{equation}\label{eqn:psiphiroughpf2}
m_\ell(2^{-k-\ell}\xi_1, 2^{-k}\eta) = \int_{\R^2} e^{-2\pi i (u2^{-k-\ell}\xi_1 + iv2^{-k}\eta)}  \widecheck{m_\ell}(u,v) d(u,v).
\end{equation}
Suppose for a moment that we can prove the decay estimate
\begin{equation}\label{eqn:kernelest}
 |\widecheck{m_\ell}(u,v)| \lesssim_N 2^{-\ell} (1+|v|)^{-2} (1+|u|)^{-N}
\end{equation}
for all $\ell>\ell_0$. Then we insert \eqref{eqn:psiphiroughpf2} into \eqref{eqn:psiphiroughpf3} and fix $u,v\in\R$. The corresponding contribution to the multilinear form has the symbol
\[ \sum_{k\in\Z} \Big(\psi_{k+\ell}(\xi_1)e^{-2\pi i u2^{-(k+\ell)}\xi_1}\Big) \Big(\psi_{k}(\xi_2)e^{-iv2^{-k}\alpha_2 \xi_2}\Big) \Big(\prod_{j=3}^{n-1} \varphi_{k}(\xi_j) e^{-iv2^{-k}\alpha_{n-1}\xi_{n-1}} \Big) \Big(\varphi_k(\xi_n) e^{-iv 2^{-k}\xi_n}\Big).
 \]
Repeating the same argument as in the $\psi$-$\psi$ case (more specifically, repeating the steps indicated from \eqref{eqn:offcoordpf3} to \eqref{eqn:offcoordpf5}), we obtain a bound that depends only logarithmically on $u,v$ and does not depend on $\ell$ (it does not depend on $\ell$, because one can shift each index $k$ appearing in \eqref{eqn:offcoordpf6} separately). The estimate \eqref{eqn:kernelest} then allows us to sum these bounds in $\ell$ and integrate in $u,v\in\R$. It now only remains to verify \eqref{eqn:kernelest}. We write
\begin{equation}\label{eqn:kernelestpf}
\widecheck{m_\ell}(u,v) = \int_{\R^2} e^{2\pi i(u \xi_1 + v\eta)} \widetilde{\mu}(\xi_1, 2^{-\ell} \eta) \widetilde{\psi}(\xi_1) \widetilde{\varphi}(\eta) d(\xi_1,\eta). 
\end{equation}
Note that on the support of the integrand we have $\xi_1>0$, $\eta<0$ and $|2^{-\ell}\eta/\xi_1|<1$ and for such $(\xi_1,\eta)$ we have
\[ \widetilde{\mu}(\xi_1, 2^{-\ell} \eta) = -2 \int_0^1 \mathbf{1}_{<0}(\alpha + 2^{-\ell}\eta/\xi_1) d\alpha = 2^{-\ell+1} \frac{\eta}{\xi_1}. \]
On the other hand, if $\eta>0$, then $\widetilde{\mu}(\xi_1, 2^{-\ell}\eta)=0$. Therefore \eqref{eqn:kernelestpf} is equal to
\[2^{-\ell+1} \Big( \int_{-\infty}^0 e^{2\pi i v\eta} \widetilde{\varphi}(\eta) \eta  d\eta\Big)\Big(\int_{\R} e^{2\pi i u\xi_1} \widetilde{\psi}(\xi_1) \xi_1^{-1} d\xi_1\Big). \]
Integrating by parts $N$ times in $\xi_1$ yields the decaying factor $(1+|u|)^{-N}$. In the $\eta$ variable we may only integrate by parts twice, which yields
\[ \int_{-\infty}^0 e^{2\pi i v\eta} \widetilde{\varphi}(\eta) \eta  d\eta = \frac{1}{(2\pi i v)^2} \Big( -1 + \int_{-\infty}^0 e^{2\pi i v \eta} (\widetilde{\varphi}''(\eta)\eta + 2\widetilde{\varphi}'(\eta)) d\eta\Big). \]
Thus we have verified \eqref{eqn:kernelest}.

\section{\texorpdfstring{An estimate in $\R^3$}{An estimate in R3}}\label{sec:twisted}
In this section we indicate how to prove \eqref{eqn:Lpboundhm} in the case $n=3$. The argument closely follows Kova\v{c}'s proof for boundedness of the twisted paraproduct \cite{K12:tp} and here we only provide the necessary changes in his argument. By a standard cone decomposition (see \cite{Thi06}) the symbol can be written as a certain superposition of symbols of the form
\[ m_1(\xi_1,\xi_2,\xi_3) = \sum_{k\in\Z} c_k \psi_k(\xi_1) \varphi_k(\xi_2) \varphi_k(\xi_3),  \]
\[ m_2(\xi_1,\xi_2,\xi_3) = \sum_{k\in\Z} c_k \varphi_k(\xi_1) \psi_k(\xi_2) \varphi_k(\xi_3),  \]
\[ m_3(\xi_1,\xi_2,\xi_3) = \sum_{k\in\Z} c_k \varphi_k(\xi_1) \varphi_k(\xi_2) \psi_k(\xi_3),  \]
where the coefficients $c_k$ satisfy $|c_k|\le 1$, $\varphi$ is a $\varphi$-type function, and $\psi$ a $\psi$-type function as in \S \ref{sec:main}. Without loss of generality we will assume that $c_k=1$ for all $k$ (see \cite[\S 1]{K12:tp}). Each of $m_1,m_2,m_3$ is treated in the same way. We focus on $m_3$. The corresponding multilinear form can be written as
\[ \Lambda_{m_3}(F_0,F_1,F_2,F_3) = \langle T_c(F_1,F_2,F_3), F_0\rangle, \]
where
\[ T_c(F_1,F_2,F_3) = \sum_{k\in\Z} (F_1*_1\widecheck{\varphi_k}) (F_2*_2\widecheck{\varphi_k}) (F_3*_3\widecheck{\psi_k}),\]
where $*_j$ denotes convolution in the $j$th variable. We now pass to a dyadic model operator for $T_c$. By $\mathcal{D}$ we denote the collection of dyadic intervals in $\R$. We define the Haar function $\psi_I$ adapted to $I\in\mathcal{D}$ by
\[ \psi_I = |I|^{-1/2} (\mathbf{1}_{I_{\mathrm{L}}} - \mathbf{1}_{I_{\mathrm{R}}}  ), \]
where $I_\mathrm{L}, I_\mathrm{R}$ denote the left and right dyadic subintervals of $I$. Also write $\varphi_I = |I|^{-1/2} \mathbf{1}_I$. Denote dyadic martingale averages and differences by
\[ \mathbb{E}_k f = \sum_{I\in\mathcal{D}, |I|=2^{-k}} \langle f, \varphi_I\rangle \varphi_I  ,\quad \Delta_k f = \sum_{I\in\mathcal{D}, |I|=2^{-k}} \langle f,\psi_I\rangle \psi_I= \mathbb{E}_{k+1} f - \mathbb{E}_k f.\]
Here $\langle f,g\rangle$ denotes $\int fg$. If $F$ is a function on $\R^3$ and $j\in\{1,2,3\}$, we write $\mathbb{E}^j_k F$, $\Delta^j_k F$ to denote application of $\mathbb{E}_k$ or $\Delta_k$ in the $j$th variable, respectively. We define a trilinear operator by
\[ T_d (F_1,F_2,F_3) = \sum_{k\in\Z} (\mathbb{E}_k^1 F_1) ( \mathbb{E}_k^2 F_2 )  ( \Delta_k^3 F_3 ). \]
Let us suppose for a moment that we can prove
\begin{equation}\label{eqn:dyadicR3}
 \|T_d(F_1,F_2,F_3)\|_{p_0'} \lesssim \prod_{j=1}^3 \|F_j\|_{p_j}
\end{equation}
for $(p_0,p_1,p_2,p_3)$ satisfying $\sum_{j=0}^3 p_j=1$ and $p_0,p_1,p_2\in (4,\infty)$, $p_3\in [2,\infty)$ and $p_0'$ the H\"older dual exponent of $p_0$.
We can then use a square-function estimate due to Jones, Seeger and Wright \cite{JSW08} comparing martingale averages with Littlewood--Paley projections to deduce the inequality
\begin{equation}\label{eqn:contR3}
\|T_c(F_1,F_2,F_3)\|_{p_0'} \lesssim \prod_{j=1}^3 \|F_j\|_{p_j}
\end{equation}
for $(p_0,p_1,p_2,p_3)$ satisfying $\sum_{j=0}^3 p_j=1$ and $p_0,p_1,p_2\in (4,\infty)$, $p_3\in [2,\infty)$. This argument is detailed in \cite[\S 6]{K12:tp} in the case of the (two-dimensional) twisted paraproduct. Up to obvious modifications this argument also applies to our case, so we omit the details.

Using multilinear interpolation and the fiber-wise Calder\'{o}n--Zygmund decomposition due to Bernicot \cite{Ber12} we can extend the range of exponents in this estimate to obtain \eqref{eqn:contR3} for $(p_0,p_1,p_2,p_3)\in (1,\infty)$ satisfying $\sum_{j=0}^3 p_j^{-1}=1$, $p_0^{-1}<2^{-2}$, $p_3^{-1}>2^{-2}$, $p_1^{-1}+p_2^{-1}>2^{-2}$, $p_2^{-1}+p_3^{-1}>2^{-1}$, and $p_1^{-1}+p_3^{-1}>2^{-1}$. This follows in the same way as in Bernicot's original argument up to obvious modifications, so we leave it out. By duality, this implies bounds for the form $\Lambda_{m_3}$. Arguing by symmetry, we also obtain corresponding bounds for $\Lambda_{m_1}$ and $\Lambda_{m_2}$. Altogether, we obtain bounds for the form $\Lambda_m$ associated with the original symbol $m$ that are valid for exponent tuples $(p_0,p_1,p_2,p_3)$ that lie in the intersection of the three regions of exponents stemming from $\Lambda_{m_1}, \Lambda_{m_2}, \Lambda_{m_3}$, respectively. This region consists of all $(p_0,p_1,p_2,p_3)\in (4,\infty)\times (1,4)^3$ so that $\sum_{j=0}^3 p_j^{-1}=1$.

It now only remains to address the validity of the estimate \eqref{eqn:dyadicR3}. This estimate follows from the techniques in \cite{K12:tp}. The details can be found in \cite{Dur14} (see \S 4.4 there). For the reader's convenience we also sketch the estimate here in the special case that $p_3=2$. It can be proven in the local $L^2$ range $p_3\in (2,\infty)$ by considering an appropriate local form as in \cite[\S 3]{K12:tp}. For this purpose it is again convenient to consider a multilinear form. We write
\[ \Lambda_d(F_0,F_1,F_2,F_3) = \langle T_d(F_1,F_2,F_3), F_0\rangle.  \]
Without loss of generality we may assume that the functions $F_j$ are non-negative. This is by splitting both the real and imaginary parts into positive and negative parts each. Expanding the martingale averages and differences and using Fubini's theorem, we see that $\Lambda_d(F_0,F_1,F_2,F_3)$ is equal to
\[\nonumber \sum_{Q} \int_{\R^4} \Big(\int_{\R} F_0(x,y,z) F_1(x',y,z) F_2(x,y',z) \psi_{I_3} (z) dz\Big) 
\Big( \int_\R F_3(x,y,z') \psi_{I_3}(z') dz'\Big) \]
\[ \varphi_{I_1}(x) \varphi_{I_1}(x') \varphi_{I_2} (y) \varphi_{I_2} (y')  d(x,y,x',y'),  \]
where the sum is over all dyadic cubes $Q=I_1\times I_2\times I_3\subset \R^3$.
By the Cauchy--Schwarz inequality this is no greater than
\begin{equation}\label{eqn:R3cauchyschwarz1}
\Big(\sum_{Q} \int_{\R^4} \Big(\int_{\R} F_0(x,y,z) F_1(x',y,z) F_2(x,y',z) \psi_{I_3} (z) dz\Big)^2 \varphi_{I_1}(x) \varphi_{I_1}(x') \varphi_{I_2} (y) \varphi_{I_2} (y')  d(x,y,x',y')\Big)^{1/2}
\end{equation}
times
\begin{equation}\label{eqn:R3cauchyschwarz2}
\Big(\sum_{Q} \int_{\R^4} \Big( \int_\R F_3(x,y,z') \psi_{I_3}(z') dz'\Big)^2 \varphi_{I_1}(x) \varphi_{I_1}(x') \varphi_{I_2} (y) \varphi_{I_2} (y')  d(x,y,x',y')\Big)^{1/2},
\end{equation}
We will treat these two factors separately. Integrating in $x',y'$, \eqref{eqn:R3cauchyschwarz2} becomes
\[ \Big(\sum_{k\in\Z} \sum_{|I_3|=2^{-k}} \int_{\R^4} F_3(x,y,z) F_3(x,y,z') \psi_{I_3}(z)\psi_{I_3}(z') d(x,y,z,z')\Big)^{1/2},
 \]
which equals $\|F_3\|_2$. On the other hand, by Theorem 6 from \cite{K12:tp} applied to the $6$-tuple $(F_0, F_0, F_1, F_1, F_2, F_2)$ we obtain that the square of \eqref{eqn:R3cauchyschwarz1} is bounded by
\[ \lesssim \|F_0\|_{p_0}^2 \|F_1\|_{p_1}^2 \|F_2\|_{p_2}^2, \]
where $\frac1{p_0}+\frac1{p_1}+\frac1{p_2} = \frac12$ and $p_0,p_1,p_2\in (4,\infty)$.

\section{Endpoint counterexample}\label{sec:counterex}

Let $n\ge 2$ and $j_0\in \{1,\dots,n\}$. Let $p_{j_0}=\infty$ and choose the other $(p_j)_j$ from $[1,\infty]$ such that $\sum_{j=0}^n p_j^{-1}=1$. We claim that the estimate \eqref{eqn:Lpbound} fails. Let $N\gg 1$ and set
\[ F_0(x) = \mathbf{1}_{(0,1)}(x_{j_0}) \prod_{j\not=j_0} \mathbf{1}_{(0,N)}(x_j),\quad F_{j_0} (x) = \mathbf{1}_{>0}(x_{j_0}), \]
\[ F_j (x) = \mathbf{1}_{(0,1)}(x_{j_0}) \prod_{\ell\not=j_0} \mathbf{1}_{(-N,N)}(x_\ell)\quad\text{for}\;j\not=0,j_0. \]
Clearly, $\|F_{j_0}\|_{\infty} = 1$ and $\|F_j\|_{p_j} \approx N^{(n-1)/p_j}$ for $j\not=j_0$. Thus, the right hand side of \eqref{eqn:Lpbound} equals $N^{n-1}$. By a straightforward computation the left hand side of \eqref{eqn:Lpbound} is comparable to $N^{n-1}\log(N)$. Letting $N\to\infty$ implies the claim.

\section{Some comments on the multiplier}\label{sec:multiplier}

In this section we are concerned with the function $\mu_n$ defined in \eqref{eqn:multiplier}. This function is the symbol of our multilinear form \eqref{eqn:mainform} and also conincides with the symbol of the Calder\'{o}n commutator \eqref{eqn:commutators}. Apart from the quadratic Fourier decay of $\mu_2$ (see Lemma \ref{lem:fourierdecay}) the only additional information which we need in \S \ref{sec:main} is the following.
%\[ \mu_n(\xi) = \int_{[0,1]^{n-1}} \mathrm{sgn}(\alpha_1\xi_1 + \cdots + \alpha_{n-1} \xi_{n-1} + \xi_n) d\alpha, \]
%where $n\ge 2$ and $\xi\in\R^n$.

\begin{lemma}\label{lem:hyperplanes}
	The function $\mu_n$ is smooth on the complement of the hyperplanes
	\begin{equation}\label{eqn:badhyperplanes}
	\{\xi\in\R^n\,:\, \alpha_1 \xi_1 + \cdots + \alpha_{n-1}\xi_{n-1}+\xi_n = 0\}
	\end{equation}
	where $\alpha\in \{0,1\}^{n-1}$.
\end{lemma}
This is a routine verification. For convenience of the reader we include some details here.
First, note the following explicit formula
\begin{align}\label{eqn:multiplierformula}
\mu_n(\xi) = \frac{c_n}{\xi_1 \ldots  \xi_{n-1}} \sum_{\alpha\in \{0,1\}^{n-1}} (-1)^{|\alpha|} g_n\Big(\xi_n+ \sum_{j=1}^{n-1} \alpha_{j}\xi_j \Big),
\end{align}
where $g_n(t)=|t|t^{n-2}$, $c_n=\frac{(-1)^{n+1}}{(n-1)!}$,  and  $\xi=(\xi_1,\ldots, \xi_n)$,  $\alpha=(\alpha_1,\ldots , \alpha_{n-1})$. This formula holds if  $\xi_1,\ldots ,\xi_{n-1}\neq 0$. 
If $k\ge 1$ of the numbers $\xi_1,\dots,\xi_{n-1}$ are equal to zero, then
\[ \mu_n(\xi) = \mu_{n-k}(\xi_*), \]
where $\xi_*$ is the vector in $\R^{n-k}$ obtained by removing the zero coordinates among $\xi_1,\dots,\xi_{n-1}$.
This is proven by induction on $n$.
Let us denote the union of the hyperplanes \eqref{eqn:badhyperplanes} by $\mathbf{B}$. Then $\R^n\setminus\mathbf{B}$ consists of a finite number of connected components each of which is a convex polytope that we call a \emph{sector}. Each sector corresponds to a map $\sigma:\{0,1\}^{n-1}\to \{-1,1\}$ where $\sigma(\alpha)$ prescribes the sign of $\xi_n+\sum_{j=1}^{n-1}\alpha_j\xi_j$.
On the closure of each such sector $\mu_n$ is given by a homogeneous rational function with no singularities in the closure of the sector. This in particular implies that $\mu_n$ is smooth on $\R^n\setminus\mathbf{B}$.

\noindent \emph{Remark.} It is natural to ask how many of the $2^{2^{n-1}}$ maps $\sigma$ correspond to (non-empty) sectors. To this end we observe that there is a bijection between sectors and linear boolean \emph{threshold functions} on the $(n-1)$-hypercube. A threshold function on the $n$-cube $\{-1,1\}^n$ is a map $f:\{-1,1\}^n\to \{-1,1\}$ that takes the form $f(x)=\mathrm{sgn}(a+\langle x,b\rangle)$ for some $a\in\R, b\in\R^n$. Threshold functions are of practical relevance in pattern recognition and have been studied extensively in the computer science literature \cite{Mur71}. In particular, it has been shown that $2^{n^2+o(n^2)}\le \#_n\le 2^{n^2}$ (see \cite{Mur71}, \cite{Zue89}), where $\#_n$ is the number of threshold functions on $\{-1,1\}^n$. In view of the separating hyperplanes theorem, the number $\#_n$ also has a simple geometric interpretation: it equals the number of subsets $A$ of vertices of the $n$-cube such that the convex hull of $A$ is disjoint from the convex hull of the vertices not in $A$.

\end{document}